\documentclass[12pt,reqno]{amsart}
\swapnumbers
\usepackage{a4,amssymb,amsthm,amscd,amsmath, verbatim,url, enumerate,mathdots, hyperref,
cleveref,mathtools}
\usepackage[pdftex,dvipsnames]{xcolor}
\usepackage{fancyhdr}
\usepackage{amsmath} 
\usepackage{hyperref}
\hypersetup{
    colorlinks=true,
    linkcolor=blue,    
    citecolor=blue,   
    urlcolor=blue
}
\makeatletter
\renewcommand{\subsection}{%
  \@startsection{subsection}{2}{\z@}%
    {-3.25ex\@plus -1ex\@minus -.2ex}%
    {1.5ex\@plus .2ex}%
    {\normalfont\bfseries}}
\makeatother
\title{Classes of Non-Invertible Simple Derivations of $k[x_1,\ldots,x_n]$}
\author{Pankaj Shukla}
\author{Sneha Singh}

\renewcommand{\deg}{\mathsf{deg}}

\renewcommand{\setminus}{\smallsetminus}

\renewcommand{\leq}{\leqslant}
\renewcommand{\geq}{\geqslant}

\swapnumbers

\newtheorem*{thm*}{Theorem}

\newtheorem{thm}{Theorem}
\numberwithin{thm}{section}
\newtheorem{prop}[thm]{Proposition}
\newtheorem*{prop*}{Proposition}
\newtheorem{cor}[thm]{Corollary}

\newtheorem{lem}[thm]{Lemma}

\theoremstyle{definition}

\newtheorem{nota}[thm]{Notation}

\newtheorem{rem}[thm]{Remark}

\newtheorem{obs}[thm]{Observation}

\newtheorem{assume}[thm]{Assumption}

\numberwithin{equation}{thm}
\renewenvironment{proof}{\par\noindent {\em Proof:}}{\hfill$\Box$\medskip}
\theoremstyle{plain}

\begin{document}

\begin{abstract}
   Let $k$ be a field of characteristic zero, and let $n\geq 3$ be an integer. For $2\leq j\leq n-1$, let $a_j\in k$, and $g_j(x_j)\in k[x_j]\setminus k.$ We give two classes of simple and non-invertible $k$-derivations of $k[x_1,\ldots,x_n].$ The first class consists of the $k$-derivations $(1 +x_1^{m_2}x_2)\partial_{x_1}+x_1^{m_1}\partial_{x_2}+(a_2 x_3 +g_2(x_2))\partial_{x_3}+\cdots+(a_{n-1}x_n + g_{n-1}(x_{n-1}))\partial_{x_n},$ where $m_1$ and $m_2$ are two positive integers such that $m_2$ does not divide $m_1$. The second class consists of $k$-derivations $(1-x_2^{\alpha} x_1^{\beta})\partial_{x_1} +x_1^m\partial_{x_2}+(a_2 x_3 +g_2(x_2))\partial_{x_3}+\cdots+(a_{n-1}x_n + g_{n-1}(x_{n-1}))\partial_{x_n},$ where $\alpha\geq 1,\beta\geq 1$ and $m\geq \beta+1$ are integers, under the assumption that $a_j=0,$ for all $2\leq j\leq n-1$ or $\deg_{x_2}g_2(x_2)=1$. 
\end{abstract}
\maketitle
\section{Introduction}
Let $k$ be a field of characteristic zero. Let $R$ be a commutative $k$-algebra. A $k$-derivation of $R$ is a $k$-linear map $d:R\to R$ that satisfies the Leibniz rule, that is, $d(ab)=ad(b)+bd(a),$ for all $a,b\in R.$ An ideal $I$ of $R$ is called a $d$-ideal of $R$ if $d(I)\subseteq I.$ A $k$-derivation $d$ is said to be simple if the ideal $\langle0\rangle$ and $R$ are only $d$-ideals of $R$. Simple derivations have several applications, such as constructing nonholonomic irreducible modules over Weyl algebras. In particular, for a second Weyl algebra $A_2$, that is, the ring of differential operators over the polynomial ring $k[x,y]$ and a simple derivation $\delta=\partial_x+f\partial_y$ of $k[x,y]$ with $f\in k[x,y]$, it has been shown in \cite{DLR06} that there exists a polynomial $g\in k[x,y]$ such that $\frac{A_2}{A_2(\delta+g)}$ is a simple nonholonomic $A_2$-module. Simple derivations also play an important role in studying Lie rings associated with derivations. It is known from \cite{Now85} that the ring $R_0$, a Lie ring with underlying additive group $R$ and Lie bracket given by $[a,b]=ad(b)-bd(a)$ for all $a,b\in R$, is simple if and only if $d$ is a simple derivation of $R$.

Let $n\geq 2$ and let $R_n:=k[x_1,\ldots,x_n]$ be the polynomial ring in indeterminates $x_1,\ldots,x_n$ over $k.$ For $n=2,$ many examples of simple $k$-derivation of $R_2$ are known; see \cite{Hav09, Now08, Kour14, SP25}. Comparatively few examples are known in higher dimensions; see \cite{Now94, AJA01, Leq08}. In this article, we mainly focus on simple $k$-derivations of $R_n$ that are non-invertible, that is, simple $k$-derivations which do not contain any unit of $R_n.$ One motivation for studying such derivations is to explore the skew polynomial rings $R[x,d]$, a non-commutative structure built from a ring $R$ equipped with a non-invertible simple derivation $d$. In \cite{Jor81}, D.A. Jordan proved that the $k$-derivation $(1-x_1x_2)\partial_{x_1}+x_1^{3}\partial_{x_2}+x_2\partial_{x_3}+\cdots+x_{n-1}\partial_{x_n}$, of $R_n,$ is simple and non-invertible, thereby providing one of the first examples of such a derivation. In \cite{MMS26}, Mishra, Mondal, and Shukla generalize Jordan's result to give a more general class of non-invertible simple derivations $(1-x_1x_2^{\alpha})\partial_{x_1}+x_1^{m}\partial_{x_2}+x_2\partial_{x_3}+\cdots+x_{n-1}\partial_{x_n},$ where $\alpha\geq 1$ and $m\geq 2$ are integers, of $R_n$. In \cite{MMS25'}, the same authors give another general class of non-invertible simple $k$-derivations $(1+x_1^{m_2}x_2)\partial_{x_1}+x_1^{m_1}\partial_{x_2}+g_2(x_2)\partial_{x_3}+\cdots+g_{n-1}(x_{n-1})\partial_{x_n},$ where $m_1$ and $m_2$ are positive integers such that $m_2$ does not divide $m_1,$ and $g_j(x_j)\in k[x_j]\setminus k$ for all $2\leq j \leq n-1,$ of $R_n.$ One more class of such derivations is known due to Kour \cite[Theorem 5.1]{Kour14}. 

In the subsequent section, we consider two classes of simple $k$-derivations of $R_2,$ namely, $(1+x_1^{m_2}x_2)\partial_{x_1}+x_1^{m_1}\partial_{x_2},$ where $m_2$ does not divide $m_1$, and $(1-x_2^{\alpha}x_1^{\beta})\partial_{x_1}+x_1^{m}\partial_{x_2},$ where $m\geq \beta+1.$ We first show that such derivations are non-invertible. In particular, we prove a stronger result for these classes of derivations of  $R_2$. For the first derivation $d_2=(1+x_1^{m_2}x_2)\partial_{x_1}+x_1^{m_1}\partial_{x_2},$ of $R_2$, we show that there is no polynomial $f\in R_2\setminus k$ such that $d_2(f)=af+g(x_2)$, where $a\in k$ and $g(x_2)\in k[x_2].$ For the derivation, $d_2=(1-x_2^{\alpha}x_1^{\beta})\partial_{x_1}+x_1^{m}\partial_{x_2}$ of $R_2,$ we prove the analogous result under some additional hypothesis. In the last section, we extend these classes of $k$-derivations of $R_2$ to $ R_n$; $ n\geq 3$. Let $a_j\in k$ and $g_j(x_j)\in k[x_j]\setminus k$, for all $2\leq j\leq n-1,$ we consider the  following two classes of extended $k$-derivations of $R_n$: $$(1 +x_1^{m_2}x_2)\partial_{x_1}+x_1^{m_1}\partial_{x_2}+(a_2 x_3 +g_2(x_2))\partial_{x_3}+\cdots+(a_{n-1}x_n + g_{n-1}(x_{n-1}))\partial_{x_n},$$ where $m_1$ and $m_2$ be two positive integers such that $m_2$ does not divide $m_1$, and the $k$-derivation $$(1-x_2^{\alpha} x_1^{\beta})\partial_{x_1} +x_1^m\partial_{x_2}+(a_2 x_3 +g_2(x_2))\partial_{x_3}+\cdots+(a_{n-1}x_n + g_{n-1}(x_{n-1}))\partial_{x_n},$$ where $\alpha\geq 1,\beta\geq 1$ and $m\geq \beta+1$ be integers. We prove that the derivations in the first class are simple and non-invertible $k$-derivations of $R_n,$ and we prove the same for the derivations in the second class under the additional assumption that $deg_{x_2}g_2(x_2)=1$ or $a_j=0,$ for all $2\leq j\leq n-1.$ This generalizes the classes known due to Mishra, Mondal, and Shukla \cite{MMS26,MMS25'}. 

\section{Two Variable Case}
Let $k$ be a field of characteristic zero. Let $x$ and $y$ be indeterminates, and $k[x,y]$ be a polynomial ring in $x$ and $y$ over $k$. Let $d$ be a $k$-derivation of $k[x,y].$ We consider the following condition for a $k$-derivation $d$ of $k[x,y].$ 

\emph{If $a\in k$, $g(x)\in k[x]$ and $f\in k[x,y]$ such that $d(f)=af+g(x)$. Then $f\in k$ and $g(x)\in k.$} 

Note that if $d$ satisfies the above condition, then $d$ is non-invertible. This condition is useful in section \ref{The General Case}.\\

In this section, we consider two classes of $k$-derivations of $k[x,y],$ namely, $y^{m_1}\partial_{x}+(1+y^{m_2}x)\partial_{y}$ and $y^m\partial_x+(1-x^{\alpha} y^{\beta})\partial_y,$ where $m_1,m_2,m,\alpha$ and $\beta$ are positive integers. We prove that for positive integers $m_1$ and $m_2$, where $m_2$ does not divide $m_1$, the derivation $y^{m_1}\partial_{x}+(1+y^{m_2}x)\partial_{y}$ satisfies the above condition (cf., Proposition~\ref{prop:main}). For positive integers $\alpha,\beta$ and $m$ with $m\geq \beta+1$, we show that the derivation $y^m\partial_x+(1-x^{\alpha} y^{\beta})\partial_y$ satisfies the above condition when $a=0$, and when $a\neq 0$ with deg$(g(x))\leq 1$(cf., Proposition~\ref{mainprop}).   

\subsection{The Derivation $y^{m_1}\partial_{x}+(1+y^{m_2}x)\partial_{y}$ of $k[x,y]$} 
\vspace{0.5 cm}
Let $m_1,m_2$ be two positive integers such that $m_2 \nmid m_1$. Let $k[x,y]$ be the polynomial ring over $k$, where $x$ and $y$ are two indeterminates. Consider the $k$-derivation $d:=y^{m_1}\partial_x+(1+y^{m_2}x)\partial_{y}$ of $k[x,y]$.

Throughout this subsection, we adhere to the following notation.

\begin{nota}
    For any polynomial $h(y)\in k[y],$ we denote the derivative of $h(y)$ with respect to $y$ by $h'(y).$ We use the usual notation $\int h(y)dy$ to denote the antiderivative of $h(y).$ 
\end{nota}

In this subsection, we prove that the $k$-derivation $d$ is non-invertible. In particular, we establish the following stronger result, that is, the derivation $d$ satisfies the condition considered above.

The proof of the following proposition is along the same lines as \cite[Lemma 3.1]{MMS25'}. Nevertheless, we write the proof for the sake of completion.

\begin{prop}\label{prop:main}
Let $d$ be the derivation of $k[x,y]$, as defined above. Let $f \in k[x,y]$ such that $d(f)=af+g(x)$, where $a \in k$ and $g(x)\in k[x]$. Then $f\in k$ and $g(x) \in k$.
\end{prop}

\begin{proof}
Let $a\in k$ and $g(x)\in k[x]$. Suppose $f(x,y)\in k[x,y]$ be such that $d(f)=af+g(x)$. Now, if $f=0$, we are done. Thus we assume that $f\neq 0$ and write $f=\sum_{i=0}^{t}f_ix^i$, where $f_i\in k[y]$, for $0\leq i\leq t$, and $f_t\neq 0$.  
	We will prove that for $t=0$, $f\in k$ and $g(x)\in k$, and for $t\geq 1$ there does not exist any such $f$.  
	
	If $t=0$ then $f=f_0(y)\in k[y]$. Now, $d(f_0(y))=af_0(y)+g(x)$ implies $(y^{m_2}x+1)f_0'=af_0(y)+g(x)$. Thus, $f_0'(y)=0$. Hence, $d(f_0)\in k$ and $f=f_0\in k$. Therefore $g(x)=-af_0\in k$.
	
	Next, we assume that $t\geq 1$ and $g(x)=b_lx^l+b_{l-1}x^{l-1}+\dots+b_1x+b_0$, where $b_i\in k$ for all $i=1,\ldots,l$, $b_l\neq 0$. From $d(f)=af+g(x)$ it follows that
	\begin{equation}\label{eqnzero}
		\sum_{i=1}^{t}if_iy^{m_1}x^{i-1}+\sum_{i=0}^{t}f'_ix^i+\sum_{i=0}^{t}y^{m_2}f'_ix^{i+1}-\sum_{i=0}^{t}af_ix^i= b_lx^l+b_{l-1}x^{l-1}+\dots+b_1x+b_0.
	\end{equation}
	Note that $f_j=0$ for $j<0$ or $j>t$. Thus, \Cref{eqnzero} becomes
	\begin{equation}\label{maineqn}
		\sum_{i=0}^{t+1}\bigg((i+1)f_{i+1}y^{m_1}+f'_{i-1}y^{m_2}+f_i'-af_i\bigg)x^i=b_lx^l+b_{l-1}x^{l-1}+\dots+b_1x+b_0.
	\end{equation}
	Thus $t+1\geq l$. Let $t+1=l$. From \Cref{maineqn} it follows that $y^{m_2}f_t'=b_l\in k^{\ast}$, which is absurd. Hence $t+1\geq l+1$, that is, $t\geq l$. 
	
	Since $t\geq l$, that is, $l\in \{t,t-1,t-2,\dots,1,0\}$, we rewrite $g(x)=b_{t}x^{t}+b_{t-1}x^{t-1}+\dots+b_1x+b_0$, where $b_i\in k $ for $0\leq i\leq t$. Similarly, rewriting \Cref{maineqn} we have
	\begin{equation}\label{maineqn2}
		\sum_{i=0}^{t+1}\bigg((i+1)f_{i+1}y^{m_1}+f'_{i-1}y^{m_2}+f_i'-af_i\bigg)x^i=b_tx^t+b_{t-1}x^{t-1}+\dots+b_1x+b_0.
	\end{equation}
	Note that the above equation also holds for $g=0$.    
	Now, comparing the coefficients of $x^i$ for $i =t+1,\ldots ,0$ in \Cref{maineqn2}, we get the following equations
	\begin{align}
		\label{eqn1}	y^{m_2}f'_t(y)&=0,\\
		\label{eqn2}	y^{m_2}f'_{t-1}(y)+f'_t(y)-af_t&=b_t,\\
		\label{eqn3}  ty^{m_1}f_t(y)+y^{m_2}f'_{t-2}(y)+f'_{t-1}(y)-af_{t-1}&=b_{t-1},\\
		\hspace{3cm} & \vdots  \nonumber
		\\
		\label{eqn4} 3y^{m_1}f_3(y)+y^{m_2}f'_{1}(y)+f'_{2}(y)-af_2&=b_2,\\
		\label{eqn5} 2y^{m_1}f_2(y)+y^{m_2}f'_{0}(y)+f'_{1}(y)-af_1&=b_1,\\
		\label{eqn6} y^{m_1}f_1(y)+f'_{0}(y)-af_0&=b_0,
	\end{align} respectively. 
	Since $t\geq 1$, we always have Equations (\ref{eqn1}) and (\ref{eqn2}), then it follows that $f_t\in k^{\ast}$ and $y^{m_2}f'_{t-1}=b_t+af_t\in k$, which implies $f_{t-1}\in k$. Thus, for $t\geq 1$, we always have that $f_t\in k^{\ast}$ and $f_{t-1}\in k$.
	Suppose $t=1$. Then from \Cref{eqn6}, we get that $y^{m_1}f_1=b_0+af_0\in k$, which gives $f_1=0$, a contradiction. Hence $t\geq 2$.

	Now, from \Cref{eqn3}, we get that
	\begin{align*}
		y^{m_2}f'_{t-2}(y)&=b_{t-1}+af_{t-1}-ty^{m_1}f_t(y),
	\end{align*}which implies that $b_{t-1}+af_{t-1}=0$. Thus
	\begin{align}\label{Eqn1}
		y^{m_2}f'_{t-2}(y)=-tf_ty^{m_1}.
	\end{align}
	
	Suppose $m_1<m_2$ then it follows from \Cref{Eqn1} that $y$ divides $f_{t}$, which is a contradiction. Hence, from now on, we consider that $m_1\geq m_2$.
	
	Since $m_2 \nmid m_1$, we have $m_1>m_2$. Then the \Cref{Eqn1} implies that
	\begin{align}\label{Eqn1'}
		f'_{t-2}(y)=-tf_ty^{m_1-m_2}.
	\end{align}
    Consider the case $t=2.$ Then \Cref{Eqn1'} becomes 
    \begin{equation}\label{eqn t=2}
       f'_{0}(y)=-2f_2y^{m_1-m_2}. 
    \end{equation} On integrating both sides of \Cref{eqn t=2} with respect to $y$, we have that
    \begin{equation}\label{eqn' t=2}
        f_0(y)=-\frac{2f_2}{m_1-m_2+1}y^{m_1-m_2+1}+c_0,
    \end{equation} where $c_0\in k.$ Now, substituting $f'_0(y)$ and $f_0(y)$ from \Cref{eqn t=2} and \Cref{eqn' t=2} into \Cref{eqn6}, respectively, we have that
    \begin{equation}\label{eqn'' t=2}
y^{m_1}f_1(y)-2f_2y^{m_1-m_2}+\frac{2af_2}{m_1-m_2+1}y^{m_1-m_2+1}-ac_0=b_0.
    \end{equation} Since $0<m_1-m_2<m_1$, the non-constant term $-2f_2y^{m_1-m_2}$ cannot be cancelled by any other non-constant terms in \Cref{eqn'' t=2}, which gives a contradiction. Hence $t\neq 2.$

    Henceforth, we consider $t\geq 3.$ On integrating \Cref{Eqn1'} with respect to $y$ both sides, we have 
	\begin{align}\label{Eqn1''}
		f_{t-2}(y)=\frac{-tf_t}{m_1-m_2+1}y^{m_1-m_2+1}+c_{t-2},
	\end{align} where $c_{t-2}\in k$.
	Next, on comparing coefficient of $x^{t-2}$ from \Cref{maineqn2}, we get that
	\begin{align}\label{eqn(t-2)}
		(t-1)y^{m_1}f_{t-1}(y)+y^{m_2}f'_{t-3}(y)+f'_{t-2}(y)-af_{t-2}=b_{t-2}.
	\end{align} Substituting $f'_{t-2}$ from \Cref{Eqn1'} and $f_{t-2}$ from \Cref{Eqn1''} in \Cref{eqn(t-2)}, we have
	\begin{align*}
		y^{m_2}f'_{t-3}(y)=b_{t-2}+ac_{t-2}-(t-1)f_{t-1}y^{m_1}-\frac{atf_t}{m_1-m_2+1}y^{m_1-m_2+1}+tf_ty^{m_1-m_2},
	\end{align*}which further implies that $b_{t-2}+ac_{t-2}=0$. Thus \Cref{eqn(t-2)} becomes
	
	\begin{align}\label{Eqn2'}
		y^{m_2}f'_{t-3}(y)=-(t-1)f_{t-1}y^{m_1}-\frac{atf_t}{m_1-m_2+1}y^{m_1-m_2+1}+tf_ty^{m_1-m_2}.
	\end{align}
	For the ease of calculation, we write 
	\begin{align}\label{Eqn2''}
		y^{m_2}f'_{t-3}(y)=A_{t-3}(y)+tf_ty^{m_1-m_2},
	\end{align} where $A_{t-3}(y)=-(t-1)f_{t-1}y^{m_1}-\frac{atf_t}{m_1-m_2+1}y^{m_1-m_2+1}$. Note that every term in $A_{t-3}(y)$ has degree strictly greater than $m_1-m_2$.
	
	Next, we claim the following for $t\geq 3$.
	
	\textbf{Claim:} Suppose $m_1>(j-1)m_2$ for some $2\leq j \leq t-1$, then we have 
	\begin{align}\label{claim-old}
		y^{m_2}f'_{t-(j+1)}=A_{t-(j+1)}(y)+(-1)^{j}tf_ty^{m_1-(j-1)m_2},
	\end{align} 
	where $A_{t-(j+1)}(y)\in k[y]$ is such that degree of each term in $A_{t-(j+1)}(y)$ is strictly greater than $m_1-(j-1)m_2$.
	
	\vspace{2mm}
		We prove the above claim by induction on $j$. \\
	\textbf{Proof of Claim:} Suppose $j=2$, that is, $m_1>m_2$ , then as in \Cref{Eqn2''} we have that $y^{m_2}f'_{t-3}(y)=A_{t-3}(y)+tf_ty^{m_1-m_2}$, where  $A_{t-3}(y)\in k[y]$ is such that degree of each term in $A_{t-3}(y)$ is strictly greater than $m_1-m_2.$ 
	Thus, the claim holds for $j=2$.

    Next, suppose for some $2\leq s-1\leq t-2$ the claim holds.  
	
	Now, assume that $m_1>(s-1)m_2$. From \Cref{maineqn2}, equating the coefficient of $x^{t-s}$, we get that

	\begin{align}\label{eqn7}
		(t-(s-1))y^{m_1}f_{t-(s-1)}(y)+y^{m_2}f'_{t-(s+1)}(y)+f'_{t-s}(y)-af_{t-s}=b_{t-s}. 
	\end{align}
	
	Note that $m_1>(s-1)m_2>(s-2)m_2$. By the induction hypothesis, we have
	
	\begin{equation}\label{eqn8-old}
		y^{m_2}f'_{t-s}(y)=A_{t-s}(y)+(-1)^{s-1}tf_ty^{m_1-(s-2)m_2},
	\end{equation}
	where $A_{t-s}(y)\in k[y]$ is such that each term in it has degree strictly greater than $m_1-(s-2)m_2$.
	Now, since $m_1>(s-1)m_2$, we can write \Cref{eqn8-old} as 
	
	\begin{equation}\label{eqn9-old}
		f'_{t-s}=B_{t-s}(y)+(-1)^{s-1}tf_ty^{m_1-(s-1)m_2},
	\end{equation}where $B_{t-s}(y)\in k[y]$ is such that $y^{m_2}B_{t-s}(y)=A_{t-s}(y),$ and each term in it has degree strictly greater than $m_1-(s-1)m_2$. Now, on integrating \Cref{eqn9-old} we have that
    \begin{equation}\label{eqn9'-old}
		f_{t-s}=\int B_{t-s}(y)dy+\frac{(-1)^{s-1}tf_t}{m_1-(s-1)m_2+1}y^{m_1-(s-1)m_2+1}+c_{t-s},
	\end{equation} where $\int B_{t-s}(y)dy\in k[y]$ is such that each term in it has degree strictly greater than $m_1-(s-1)m_2+1$ and $c_{t-s}\in k.$
    Substituting the value of $f'_{t-s}(y)$ and $f_{t-s}(y)$ from \Cref{eqn9-old} and \Cref{eqn9'-old} respectively in \Cref{eqn7}, we have that 
	
	\begin{align}\label{eqn10-old}
		y^{m_2}f'_{t-(s+1)}(y)&=b_{t-s}+ac_{t-s}-(t-(s-1))y^{m_1}f_{t-(s-1)} \nonumber\\
        &-\bigg(B_{t-s}(y)+(-1)^{s-1}tf_ty^{m_1-(s-1)m_2}\bigg)+a\int B_{t-s}(y)dy \nonumber\\
        &+\frac{(-1)^{s-1}tf_ta}{m_1-(s-1)m_2+1}y^{m_1-(s-1)m_2+1}.
	\end{align}
	Since $m_1>(s-1)m_2$, it follows that $b_{t-s}+ac_{t-s}=0$. Thus, every term in \Cref{eqn10-old} other than $(-1)^{s}tf_ty^{m_1-(s-1)m_2}$ has degree strictly greater than $m_1-(s-1)m_2.$ So, we write \Cref{eqn10-old} as
	\begin{align}\label{eqn11-old}
		y^{m_2}f'_{t-(s+1)}=A_{t-(s+1)}(y)+(-1)^{s}tf_ty^{m_1-(s-1)m_2},
	\end{align}
	where $A_{t-(s+1)}(y)\in k[y]$ and each term in it has degree strictly greater than $m_1-(s-1)m_2$. 
	Hence, the claim follows for $s$.\\
    Thus, if $jm_2>m_1>(j-1)m_2$ for some $2\leq j\leq t-1$, then from \Cref{claim-old} it follows that $y$ divides $tf_t\in k^{\ast}$, which  is a contradiction. Hence, we get that $m_1>(t-1)m_2$. Now, from \Cref{eqn6} we have 
	\begin{align}\label{eqn12}
		y^{m_1}f_1(y)+f'_0(y)-af_0=b_0. 
	\end{align}
	Also, as $m_1>(t-1)m_2>(t-2)m_2$, from \Cref{claim-old} we have 
	\begin{equation}\label{eqn13}
		f'_0(y)=B_0(y)+(-1)^{t-1}tf_ty^{m_1-(t-1)m_2},
	\end{equation}
	where $B_0(y)\in k[y]$ is such that $y^{m_2}B_0(y)=A_0(y)$, and each term in it has degree strictly greater than $m_1-(t-1)m_2$. Thus, on integrating \Cref{eqn13} both sides with respect to $y$, we have
	\begin{equation}\label{eqn14}
		f_0(y)=\int B_0(y)dy+\frac{(-1)^{t-1}tf_t}{m_1-(t-1)m_2+1}y^{m_1-(t-1)m_2+1}+c_0,
	\end{equation} where $\int B_0(y)dy\in k[y]$ is such that each term in it has degree strictly greater than $m_1-(t-1)m_2+1$ and $c_0\in k$. Thus,
	substituting expressions of $f'_0(y)$ and $f_0(y)$ from \Cref{eqn13} and \Cref{eqn14} respectively in \Cref{eqn12} it follows that $b_0+ac_0=0$. Furthermore, we get that 
	\begin{equation}
		y^{m_1}f_1(y)=A(y)+(-1)^{t}tf_ty^{m_1-(t-1)m_2},
	\end{equation} where $A(y)\in k[y]$ is such that each term in it has degree strictly greater than $m_1-(t-1)m_2$.
	Then $y$ divides $tf_t\in k^{\ast}$, which is a contradiction. Hence $t\ngeq 1$. Thus, we finally get that $t=0$ is the only possible case. This completes the proof.  
\end{proof}

\subsection{The Derivation $y^m\partial_x+(1-x^{\alpha} y^{\beta})\partial_y$}

Consider the derivation $d=y^m\partial_x+(1-x^{\alpha} y^{\beta})\partial_y$ of $k[x,y]$ with $\alpha\geq 1, \beta\geq 1$ and $m\geq \beta+1$. In this subsection, we prove that if there exists some polynomial $f\in k[x,y]$ such that $d(f)=g(x)$ for some $g(x)\in k[x]$, then $f\in k$ and $g(x)=0$(c.f. proposition~\ref{pc=0}), and if $f\in k[x,y]$ is such that $d(f)=af+g(x)$, where $a\in k^\ast$ and $g(x)\in k[x]$ with $deg_x g(x)\leq 1$, then $f\in k$ and $g(x)\in k$(c.f. proposition~\ref{mainprop}).

\noindent In this subsection, we adhere to the following notation.
\begin{nota}\label{nota:1}
   Let $\ell>0$ and $f=\sum_{i=0}^{\ell}f_i(x)y^i\in k[x,y]$ be such that $f_i(x)\in k[x]$ for $0\le i \leq\ell$, $f_\ell(x) \neq 0$ and $d(f)=af+g(x)$. Set $f_i(x)=0$ for all integers $i>\ell$ and $i<0$. Also, for any $h(x)\in k[x]$, we denote by $h'(x)$ the usual derivative of $h(x)$ with respect to $x$.
\end{nota}
\noindent Note that $\ell=0$ case will be handled later in Proposition \ref{pc=0} and Proposition \ref{mainprop}.
\begin{obs}\label{obs:1}
    Let $\ell>0$ be an integer and consider the expression $\ell-(m-(\beta-1))(j-1)$. For $j=1,\; \ell-(m-(\beta-1))(j-1)=\ell>0 $. Thus there exists one positive integer for which the expression has positive value. As $m,\beta$ and $\ell$ are positive integers, thus for $j>>0$, $\ell-(m-(\beta-1))(j-1)\leq0$. We denote by $j_0$ the largest positive integer such that $\ell-(m-(\beta-1))(j-1)>0$.
\end{obs}

\begin{lem}\label{deglem}
Let $g(x),f(x,y),\ell$ and $f_i(x)$ be as in Notation~\ref{nota:1}. Let $j_0$ be the greatest integer such that $\ell-(m-(\beta-1))(j-1)>0$. Then $\ell=(m-(\beta-1))j_0$, and for every $j$ with $1\leq j\leq j_0$ and every $s$ with $1\leq s\leq m-\beta$, $$\deg_x f_{\ell-(m-(\beta-1))j}(x)=j(\alpha+1)$$
and
$$\deg_x f_{\ell-(m-(\beta-1))j+s}(x)\leq (j-1)(\alpha+1).$$
\end{lem}
\begin{proof}
Since $f=\sum_{i=0}^{\ell}f_i(x)y^i$, where $f_i(x)\in k[x]$ with $0\leq i \leq \ell$ and $f_\ell(x)\ne 0$, and from $d(f)=af+g(x)$, we have
\begin{equation}
(y^m\partial_x+(1-x^\alpha y^\beta)\partial_y)
\sum_{i=0}^{\ell}f_i(x)y^i=a\sum_{i=0}^{\ell}f_i(x)y^i +g(x).
\end{equation}
Since $f_i(x)=0$ for $i>\ell$ or $i<0$, we get
\begin{equation}\label{eq1:2}
\sum_{i=0}^{\ell+m}
\left(f'_{i-m}(x)+(i+1)f_{i+1}(x)-(i-(\beta-1))x^\alpha f_{i-(\beta-1)}(x)-af_i(x)\right)y^i
=g(x).
\end{equation}
From the above equation, for $i\geq1$, we have
\begin{equation}\label{eq1:3}
f'_{i-m}(x)=(i-(\beta-1))x^\alpha f_{i-(\beta-1)}(x)-(i+1)f_{i+1}(x)+af_i(x).
\end{equation}
\noindent  We first determine $f_\ell, \ f_{\ell-1},\ldots, \ f_{\ell-m+\beta}$.  
Substituting $i=\ell+m,\ell+m-1,\dots,\ell+\beta$ in \Cref{eq1:3}, we have that
$$f'_\ell(x)=f'_{\ell-1}(x)=\cdots=f'_{\ell-m+\beta}(x)=0.$$
Hence,
\begin{equation}\label{eq1:4}
f_\ell(x),f_{\ell-1}(x),\dots,f_{\ell-m+\beta}(x)\in k.
\end{equation}
Without loss of generality, we may assume that $f_\ell(x)=1$.  
For $1\leq s\leq m-\beta$, let $f_{\ell-s}(x):=\lambda_s$ with $\lambda_s\in k$.

\medskip 

\noindent Substituting $i=\ell+\beta-1$ in \Cref{eq1:3}, we have that
\begin{equation}\label{eq1:5}
    f'_{\ell-m+(\beta-1)}(x)=\ell x^\alpha f_\ell(x)-(\ell+\beta)f_{\ell+\beta}(x)+af_{\ell+\beta-1}(x).
\end{equation}
Recall that $f_{\ell}=1$, and $f_i(x)=0$ for all $i>\ell$ and $i<0$. Thus $f_{\ell+\beta}=0$ and $f_{\ell+\beta-1}\in k$, since $\beta \geq 1$. On integrating \Cref{eq1:5} with respect to $x$, we have that
\begin{equation}\label{eq1:6}
f_{\ell-m+(\beta-1)}(x)=\frac{\ell x^{\alpha+1}}{\alpha+1}+af_{\ell+\beta-1}x+c_{\ell-m+(\beta-1)},
\end{equation} where $c_{\ell-m+(\beta-1)}\in k.$
 Thus, $\deg_x f_{\ell-m+(\beta-1)}(x)=\alpha+1$.

\medskip

\noindent Since $j_0\geq 1$. We first show that the lemma holds for $j_0=1$.
Assume that $j_0=1$. Then $\ell\leq (m-(\beta-1))$. If $\ell<(m-(\beta-1))$, then $f_{\ell-m+(\beta-1)}(x)=0$, 
which is a contradiction to the fact that $\deg_x f_{\ell-m+(\beta-1)}(x)=\alpha+1$. 
Thus $\ell=m-(\beta-1)$. Hence, from Equations (\ref{eq1:4}) and (\ref{eq1:6}), we have that 
$\deg_x f_s(x)\leq0$ for all $1\leq s\leq m-\beta$ and 
$\deg_x f_0(x)=\alpha+1$. Hence the lemma holds for $j_0=1$.

\medskip

\noindent Henceforth, we assume that $j_0\geq2$. Now, we show that for every 
$j$ with $1\leq j\leq j_0$, and every $s$ with $1\leq s\leq m-\beta$,
\[
\deg_x f_{\ell-(m-(\beta-1))j}(x)=j(\alpha+1)
\]
and
\[
\deg_x f_{\ell-(m-(\beta-1))j+s}(x)\leq (j-1)(\alpha+1).
\]
\noindent We will use induction on $j$ to establish the above degree inequalities.\\
From the preceding observations, we have
\[
f_\ell=1,\quad
f_{\ell-1}=\lambda_1,\quad
\dots,\quad
f_{\ell-m+\beta}=\lambda_{m-\beta},
\]
where $\lambda_1,\dots,\lambda_{m-\beta}\in k$ and $\deg_x f_{\ell-m+(\beta-1)}(x)=\alpha+1$. Hence the statement holds for $j=1$.\\
Now, let $t$ be a positive integer such that $1\leq t\leq j_0-1$. Assume that, for all $j$ with $1\leq j\leq t$,
$\deg_x f_{\ell-(m-(\beta-1))j}(x)=j(\alpha+1)$
and
$\deg_x f_{\ell-(m-(\beta-1))j+s}(x)\leq (j-1)(\alpha+1)$
for all $1\leq s\leq m-\beta$.\\
We show that $$\deg_x f_{\ell-(m-(\beta-1))(t+1)}(x)=(t+1)(\alpha+1)$$ and $$\deg_x f_{\ell-(m-(\beta-1))(t+1)+s}(x)\leq t(\alpha+1),$$ for all $1\leq s\leq m-\beta$.

\medskip
\noindent We first show that $\deg_x f_{\ell-(m-(\beta-1))(t+1)}(x)=(t+1)(\alpha+1)$.

\medskip
\noindent Substituting $i=\ell-t(m-(\beta-1))+\beta-1$ in \Cref{eq1:3}, we have that
\begin{align}\label{eq1:7}
f'_{\ell-(t+1)(m-(\beta-1))}(x)&=(\ell-t(m-(\beta-1)))x^\alpha f_{\ell-t(m-(\beta-1))}(x)\notag\\
&-(\ell-t(m-(\beta-1))+\beta)f_{\ell-t(m-(\beta-1))+\beta}(x)\notag \\
&+af_{\ell-t(m-(\beta-1))+\beta-1}(x).
\end{align}
By the induction hypothesis, $\deg_x f_{\ell-t(m-(\beta-1))}(x)=t(\alpha+1)$. We now determine the degrees of the terms $f_{\ell-t(m-(\beta-1))+\beta}(x)$ and $f_{\ell-t(m-(\beta-1))+\beta-1}(x)$ in \Cref{eq1:7}. Recall that $m\geq \beta+1.$
Before considering the following cases, note that $f_{\ell}=1$ and $f_i=0$ for all $i>\ell.$ Thus, for integers $q\leq 0$ and $0\leq r\leq m-\beta,$ we have that $f_{\ell-q(m-(\beta-1))+r}(x)\in k.$ We will use this fact throughout the proof. In the following cases, we therefore only consider the terms for which $q>0.$

\medskip
\noindent\textbf{Case 1.} When $\beta+1\leq m \leq 2\beta-3$, we have that $2\leq m-\beta+1\leq \beta-1.$ Then by the division algorithm there exists $t_1,t_2\in \mathbb{N}$ and $0 \leq r_1,r_2\leq m-\beta$ such that  $\beta=t_1(m-\beta+1)+r_1$ and $\beta-1=t_2(m-\beta+1)+r_2.$ Consequently, $f_{\ell-t(m-(\beta-1))+\beta}=f_{\ell-(t-t_1)(m-(\beta-1))+r_1}$ and $f_{\ell-t(m-(\beta-1))+\beta-1}(x)=f_{\ell-(t-t_2)(m-(\beta-1))+r_2}(x).$   Now, for $j=1,2,$ the induction hypothesis gives, $$\deg_x f_{\ell-(t-t_j)(m-(\beta-1))+r_j} \leq (t-t_j-1)(\alpha+1)$$ for $1\leq r_j\leq m-\beta,$ and $$\deg_x f_{\ell-(t-t_j)(m-(\beta-1))}(x)=(t-t_j)(\alpha+1).$$\\
\textbf{Case 2.} When $m=2 \beta -2$. We have that $\deg_x f_{\ell-t(m-(\beta-1))+\beta}$ $=\deg_x f_{\ell-(t-1)(m-(\beta-1))+1},$ and $\deg_xf_{\ell-t(m-(\beta-1))+\beta-1}(x)=\deg_x$ $f_{\ell-(t-1)(m-(\beta-1))}(x).$ We have from the induction hypothesis, $\deg_x f_{\ell-(t-1)(m-(\beta-1))+1}\leq (t-2)(\alpha+1),$ and $\deg_xf_{\ell-t(m-(\beta-1))+\beta-1}(x)=(t-1)(\alpha+1).$ 

\medskip
\noindent
\textbf{Case 3.} When $m=2\beta -1$. By the induction hypothesis, $\deg_x f_{\ell-t(m-(\beta-1))+\beta}$ $=\deg_x f_{\ell-(t-1)(m-(\beta-1))}= (t-1)(\alpha+1)$ and $\deg_xf_{\ell-t(m-(\beta-1))+\beta-1}(x)\leq (t-1)(\alpha+1).$

\medskip
\noindent
\textbf{Case 4.} When $m \geq 2\beta$. By the induction hypothesis, we have $\deg_x f_{\ell-t(m-(\beta-1))+\beta}$ $\leq (t-1)(\alpha+1)$, $\deg_xf_{\ell-t(m-(\beta-1))+\beta-1}(x)\leq (t-1)(\alpha+1),$ if $\beta>1$, and $\deg_x$ $f_{\ell-t(m-(\beta-1))+\beta-1}(x)=\deg_xf_{\ell-tm}(x)=t(\alpha+1),$ if $\beta=1.$

\medskip
\noindent Thus,  $\deg_x(x^\alpha f_{\ell-t(m-(\beta-1))}(x))$ is strictly greater than $\deg_xf_{\ell-t(m-(\beta-1))+\beta}(x)$ and $\deg_xf_{\ell-t(m-(\beta-1))+\beta-1}(x)$ in \Cref{eq1:7}. Then, integrating \Cref{eq1:7} with respect to $x$, we have $\deg_x f_{\ell-(m-(\beta-1))(t+1)}(x)=(t+1)(\alpha+1)$.

\medskip
\noindent Next, we show that $\deg_x f_{\ell-(m-(\beta-1))(t+1)+s}(x)\leq t(\alpha+1),$ for all $1\leq s\leq m-\beta.$ 
Substituting $i=\ell-t(m-(\beta-1))+s+\beta-1$, where $1\leq s\leq m-\beta$, in \Cref{eq1:3}, we get
\begin{align}\label{eq1:8}
f'_{\ell-(t+1)(m-(\beta-1))+s}(x)&=(\ell-t(m-(\beta-1))+s)x^\alpha f_{\ell-t(m-(\beta-1))+s}(x)\notag \\
&-(\ell-t(m-(\beta-1))+s+\beta)f_{\ell-t(m-(\beta-1))+s+\beta}(x) \notag \\
&+af_{\ell-t(m-(\beta-1))+s+\beta-1}(x).
\end{align}
\noindent
By the induction hypothesis, $\deg_x f_{\ell-t(m-(\beta-1))+s}(x)\leq(t-1)(\alpha+1)$. We consider the following cases to determine the degrees of the terms $f_{\ell-t(m-(\beta-1))+s+\beta}(x)$ and $f_{\ell-t(m-(\beta-1))+s+\beta-1}(x)$ in \Cref{eq1:8}. Recall that  for integers $q\leq 0$ and $0\leq r\leq m-\beta,$ we have that $f_{\ell-q(m-(\beta-1))+r}(x)\in k.$

\medskip
\noindent
\textbf{Case 1.}  Suppose $1\leq s\leq m-2\beta$. Then from the induction hypothesis, $\deg_x f_{\ell-t(m-(\beta-1))+s+\beta}(x)$, $\deg_x f_{\ell-t(m-(\beta-1))+s+\beta-1}(x)\leq (t-1)(\alpha+1)$.

\medskip
\noindent
\textbf{Case 2.} Suppose $s=m-2\beta+1$. Then $f_{\ell-t(m-(\beta-1))+s+\beta}(x)=f_{\ell-(t-1)(m-(\beta-1))}(x)$. Thus, from the induction hypothesis, $\deg_x f_{\ell-t(m-(\beta-1))+s+\beta}= (t-1)(\alpha +1)$ and $\deg_x f_{\ell-t(m-(\beta-1))+s+\beta-1}(x)=\deg_x f_{\ell-t(m-(\beta-1))+m-\beta}(x)\leq(t-1)(\alpha+1)$.

\medskip
\noindent
\textbf{Case 3.} Suppose $s=m-2\beta+2$. Then from the induction hypothesis, $\deg_x f_{\ell-t(m-(\beta-1))+s+\beta}(x)$ $=\deg_x f_{\ell-(t-1)(m-(\beta-1))+1}(x)\leq (t-2)(\alpha+1)$ and \\ $\deg_x f_{\ell-t(m-(\beta-1))+s+\beta-1}(x)=\deg_x f_{\ell-(t-1)(m-(\beta-1))}(x)=(t-1)(\alpha+1).$ 

\medskip
\noindent
\textbf{Case 4.} Suppose $m-2\beta+3\leq s\leq m-\beta$. Then by the division algorithm, there exist $t_1,t_2\in \mathbb{N}$ and $0\leq r_1 ,r_2\leq m-\beta$ such that $s+\beta=t_1(m-(\beta-1))+r_1$ and $s+\beta-1=t_2(m-(\beta-1))+r_2.$ Hence, $f_{\ell-t(m-(\beta-1))+s+\beta}=f_{\ell-(t-t_1)(m-(\beta-1))+r_1}$ and $f_{\ell-t(m-(\beta-1))+s+\beta-1}=f_{\ell-(t-t_2)(m-(\beta-1))+r_2}.$ Now, for $j=1,2$, the induction hypothesis gives, $$\deg_x f_{\ell-(t-t_j)(m-(\beta-1))+r_j} \leq (t-t_j-1)(\alpha+1)$$ for $1\leq r_j\leq m-\beta,$ and $$\deg_x f_{\ell-(t-t_j)(m-(\beta-1))+r_j}(x)=(t-t_j)(\alpha+1)$$ for $r_j=0.$\\
Consequently, the degree of the right-hand side of \Cref{eq1:8} is at most $t(\alpha+1)$. Therefore, $\deg_x f_{\ell-(t+1)(m-(\beta-1))+s}\leq t(\alpha+1).$

\medskip

\noindent Finally, we show that $\ell=(m-(\beta-1))j_0$. By the choice of $j_0$, it follows that $\ell-(m-(\beta-1))j_0\le0$. 
Suppose that $\ell-(m-(\beta-1))j_0<0$. Then $f_{\ell-(m-(\beta-1))j_0}(x)=0$. But we have $\deg_x f_{\ell-(m-(\beta-1))j_0}(x)=j_0(\alpha+1)>0.$ Thus, we get a contradiction. Hence, $\ell=(m-(\beta-1))j_0$.
\medskip

\noindent This completes the proof.
\end{proof}

\begin{assume}
For Lemma \ref{fval} and Proposition \ref{pc=0}, we consider the case $a=0$, that is, 
$d(f)=g(x)$, where $g(x)\in k[x]$. Note that $f$, $\ell$, $g(x)$ and $f_i(x)$ remain the same as in Notation \ref{nota:1} and Observation \ref{obs:1}.
\end{assume}
\begin{lem}\label{fval}
Let $g(x),f(x,y)$, $\ell$, and $f_i(x)$ be as in Notation~\ref{nota:1}. Then for $j \in \mathbb{Z}_{\geq 0}$, $f_{\beta j + i}(x) = 0\;\text{for } 2 \leq i \leq \beta$, and $f_{\beta j + 1}(x) = \frac{x^{j\alpha} g(x)}{\beta j + 1}$,
whenever $1 \leq \beta j + i \leq m$ for $1 \leq i \leq \beta$.
\end{lem}
\begin{proof}
Comparing the degree zero term of $y$ in \Cref{eq1:2}, we obtain
\begin{equation}\label{eq2:1}
f_1(x)+f'_{-m}(x)=g(x) .
\end{equation}
Hence, $$f_1(x)=g(x).$$
Further, for $i=1,\dots ,\beta-1$, from \Cref{eq1:3} we get
\begin{equation}\label{eq2:3}
f'_{i-m}(x)=(i-(\beta-1))x^\alpha f_{i-(\beta-1)}(x)-(i+1)f_{i+1}(x).
\end{equation}
For $i=1,\dots,\beta-1$, we have $(i-(\beta-1))x^\alpha f_{i-(\beta-1)}(x)=0,$ and $f'_{i-m}=0$ since $m\geq \beta+1,$ in \Cref{eq2:3}. Hence, we have that $$f_2,f_3,\dots,f_\beta=0.$$
Using this relation inductively, we have that
$$f_{\beta j + i}(x) = 0, \qquad \text{for } 2 \leq i \leq \beta,$$
and
$$f_{\beta j+1}(x) = \frac{x^{j\alpha}g(x)}{\beta j+1},$$
whenever $1\leq\beta j+i\leq m$ with $1\leq i\leq\beta$ and $j\in\mathbb{Z}_{\geq0}$.
\end{proof}

The following proof follows the proof of Proposition 2.4 in \cite{MMS26}. For the sake of completeness, we include it here.
\begin{prop}\label{pc=0}
Let $\alpha\geq 1$, $\beta\geq1$ and $m\geq \beta+1$ be integers. Let $R_2=k[x,y]$ with the $k$-derivation $d=y^m\partial_x+(1-x^\alpha y^\beta)\partial_y.$ Let $f\in R_2$ be such that $d(f)=g(x)$ for some $g(x)\in k[x]$. Then $f\in k$ and $g(x)=0$.
\end{prop}
\begin{proof}
Let $f\in R_2\setminus k$ be such that $d(f)=g(x)$ for some $g(x)\in k[x]$.\\
As $f\in R_2$, we can write $f=\sum_{i=0}^{\ell} f_i(x)y^i$, where $f_i(x)\in k[x] \; \text{for}\; 0\leq i \leq \ell $ and $f_\ell(x)\neq0$.\\
Suppose that $\ell=0$. Then $d(f)=d(f_0(x))=y^m f_0'(x)=g(x)$. Thus $f_0'(x)=0$, which implies $f=f_0(x)\in k$, a contradiction. Now, suppose  $g(x)=0$, then $d(f)=0\implies f\in k$ (as from Lemma \ref{lem simple}, we know that the derivation $d$ is simple). But that is again a contradiction.\\
Henceforth, we assume that $\ell>0$ and $g(x)\ne 0$. Now we are in the same setup as Lemmas \ref{deglem} and \ref{fval}, and so we will use facts and results from those lemmas.\\
As in the proof of Lemma \ref{deglem}, we have
\[
f_\ell(x),\, f_{\ell-1}(x),\, f_{\ell-2}(x),\,\dots,\,f_{\ell-m+\beta}(x)\in k.
\]
Further, without loss of generality, we assume that $f_\ell(x)=1$.\\
For $1\leq s\leq m-\beta$, let us denote $f_{\ell-s}(x)$ by $\lambda_s$.\\
Considering $j_0$ as in Lemma \ref{deglem}, we have $\ell=j_0(m-(\beta-1))$ and 
\begin{equation}\label{eq3.1}
    \deg_x f_{m-(\beta-1)}(x)=\deg_x f_{\ell-(j_0-1)(m-(\beta-1))}(x)=(j_0-1)(\alpha+1).
\end{equation}
From Lemma \ref{fval} we have, If $m-(\beta-1)$ has the form $j\beta +i$ for some $2\leq i \leq \beta$, where $j\in \mathbb{Z}_{\geq 0}$, that is, $m=(j+1)\beta+(i-1)$, then $f_{m-(\beta-1)}(x)=0$. But from \Cref{eq3.1}, we have $\deg_x f_{m-(\beta-1)}(x)=(j_0-1)(\alpha+1)$. Thus, we have a contradiction for $m=(j+1)\beta+r$, where $1\leq r \leq \beta-1$ and $j\in \mathbb{Z}_{\geq 0}$.\\
Again, from Lemma \ref{fval}, if $m-(\beta-1)$ has the form $j\beta +1$ for some $j\in \mathbb{N}$, that is, $m=(j+1)\beta$, then 
\begin{equation}\label{eq3.2}
    f_{m-(\beta-1)}(x)=\frac{x^{j\alpha} g(x)}{m-(\beta-1)}.
\end{equation}
Note that $j\neq 0$, otherwise $m=\beta$, which is not possible.\\
Thus, from Equations (\ref{eq3.1}) and (\ref{eq3.2}), we conclude that $j_0>1$.\\
Henceforth, we assume that $m-(\beta-1)=j\beta+1,$ that is, $m=(j+1)\beta$, where $j\in \mathbb{N}$.
From \Cref{eq1:3}, we have
\begin{equation}\label{eq3.3}
f'_{i-m}(x)=(i-(\beta-1))x^\alpha f_{i-(\beta-1)}(x)-(i+1)f_{i+1}(x),
\end{equation}
for all $i\geq1$.

\medskip
\noindent
\textbf{(I)} First we find the values of $f_{\ell-j(m-(\beta-1))}(x)$ and $f_{\ell-j(m-(\beta-1))-\beta}(x)$ for all $1\leq j\leq j_0-1$.\\
For $1\leq j\leq j_0-1$, substituting $i=\ell-(j-1)(m-(\beta-1))+\beta-1$ and $i=\ell-(j-1)(m-(\beta-1))-1$ in \Cref{eq3.3}, we obtain
\begin{align}\label{eq3.4}
f'_{\ell-j(m-(\beta-1))}(x) &=\bigl(\ell-(j-1)(m-(\beta-1))\bigr)x^\alpha f_{\ell-(j-1)(m-(\beta-1))}(x)\notag\\
&\qquad -\bigl(\ell-(j-1)(m-(\beta-1))+\beta \bigr)f_{\ell-(j-1)(m-(\beta-1))+\beta}(x),
\end{align}
and
\begin{align}\label{eq3.5}
f'_{\ell-j(m-(\beta-1))-\beta}(x) &= \bigl(\ell-(j-1)(m-(\beta-1))-\beta\bigr)x^\alpha f_{\ell-(j-1)(m-(\beta-1))-\beta}(x)\notag \\
& \qquad-\bigl(\ell-(j-1)(m-(\beta-1))\bigr)f_{\ell-(j-1)(m-(\beta-1))}(x),
\end{align}
respectively.\\
\noindent\textbf{Claim.} We claim that for all $j$, $1\leq j\leq j_0-1$,
\begin{equation}\label{eq3.6}
f_{\ell-j(m-(\beta-1))}(x)=A_{\ell-j(m-(\beta-1))}x^{j(\alpha+1)}+O_{\ell-j(m-(\beta-1))}(x)
\end{equation}
and
\begin{align}\label{eq3.7}
    f_{\ell-j(m-(\beta-1))-\beta}(x)&=A_{\ell-j(m-(\beta-1))-\beta}\;x^{j(\alpha+1)}+B_{\ell-j(m-(\beta-1))-\beta}\;x^{(j-1)(\alpha+1)+1}\notag \\
    &\hspace{5cm}+O_{\ell-j(m-(\beta-1))-\beta}(x),
\end{align}
where $A_{\ell-j(m-(\beta-1))}\in \mathbb{Q}_{>0}$, $B_{\ell-j(m-(\beta-1))-\beta}\in \mathbb{Q}_{<0}$, $A_{\ell-j(m-(\beta-1))-\beta}\in k$, \\
$O_{\ell-j(m-(\beta-1))}(x)\in k[x]$, $\deg_x O_{\ell-j(m-(\beta-1))}(x)<j(\alpha+1)$, and $O_{\ell-j(m-(\beta-1))-\beta}(x)\in k[x],\deg_x O_{\ell-j(m-(\beta-1))-\beta}(x)<(j-1)(\alpha+1)+1$.

\medskip
\noindent
We prove the above claim by induction on $j$.\\
Note that $f_{\ell-\beta}(x)=\lambda_\beta$ (since $ f_{\ell-s}=\lambda_s$ for $1\leq s \leq m-\beta$, and we have $m\geq 2 \beta$, which implies that $\beta\leq m-\beta$), $f_{\ell}(x)=1$ and $f_{\ell+\beta}(x)=0$.
Substituting $j=1$ in \Cref{eq3.4}, and then integrating with respect to $x$, we get
\[
f_{\ell-m+(\beta-1)}(x)=\frac{x^{\alpha+1}\ell}{\alpha+1}+c_{\ell-m+(\beta-1)},
\]
where $c_{\ell-m+(\beta-1)}\in k$.\\
Similarly, substituting $j=1$ in \Cref{eq3.5} and integrating with respect to $x$, we get
\[
f_{\ell-m-1}(x)=\frac{x^{\alpha+1}(\ell-\beta)}{\alpha+1}\lambda_\beta-\ell x+c_{\ell-m-1},
\]
where $c_{\ell-m+(\beta-1)},c_{\ell-m-1}\in k$.\\ \noindent
Note that $A_{\ell-m+(\beta-1)}=\frac{\ell}{\alpha+1}\in \mathbb{Q}_{>0},$ $B_{\ell-m-1}=-\ell\in\mathbb{Q}_{<0},$ $A_{\ell-m-1}=\frac{\ell-\beta}{\alpha+1}\lambda_\beta \in k$, $O_{\ell-m+(\beta-1)}(x)=c_{\ell-m+(\beta-1)}\in k[x]$ with $\deg_x O_{\ell-m+(\beta-1)}<(\alpha+1),$ and $O_{\ell-m-1}(x)=c_{\ell-m-1}\in k[x]$ with $\deg_x O_{\ell-m-1}<1$.
Hence, the claim holds for $j=1$.
\\ \noindent Next, as induction hypothesis we assume that the above claim holds for $(j-1)$ where $2\le j\le j_0-1$. We now show that the claim holds for $j$.
\\ \noindent Substituting the value of $f_{\ell-(j-1)(m-(\beta-1))}(x)$ from \Cref{eq3.6} in \Cref{eq3.4}, we get 
\begin{align*}
f'_{\ell-j(m-(\beta-1))}(x)
&=(\ell-(j-1)(m-(\beta-1)))x^\alpha
\left[
\begin{aligned}
&A_{\ell-(j-1)(m-(\beta-1))}
x^{(j-1)(\alpha+1)}\\
&\qquad+O_{\ell-(j-1)(m-(\beta-1))}(x)
\end{aligned}
\right]\\
&-(\ell-(j-1)(m-(\beta-1))+\beta)
f_{\ell-(j-1)(m-(\beta-1))+\beta}(x).
\end{align*}
Integrating with respect to $x$, we obtain
\begin{align*} 
f_{\ell-j(m-(\beta-1))}(x)
&=\frac{(\ell-(j-1)(m-(\beta-1))A_{\ell-(j-1)(m-(\beta-1))}}{j(\alpha+1)}x^{j(\alpha+1)}\\
& \hspace{5cm}+O_{\ell-j(m-(\beta-1))}(x),
\end{align*}
for some $O_{\ell-j(m-(\beta-1))}(x)\in k[x]$. By the induction hypothesis, $ A_{\ell-(j-1)(m-(\beta-1))}\in\mathbb{Q}_{>0},\deg_x O_{\ell-(j-1)(m-(\beta-1))}(x)<(j-1)(\alpha+1)$. Also from Lemma \ref{deglem}, we have $\deg_x f_{\ell-(j-1)(m-(\beta-1))+\beta}(x)\le (j-2)(\alpha+1)$.\\ \noindent
Hence it follows that
\[
f_{\ell-j(m-(\beta-1))}(x)
=A_{\ell-j(m-(\beta-1))}x^{j(\alpha+1)}+O_{\ell-j(m-(\beta-1))}(x),
\]
where
\[
A_{\ell-j(m-(\beta-1))}
=\frac{(\ell-(j-1)(m-(\beta-1)))A_{\ell-(j-1)(m-(\beta-1))}}{j(\alpha+1)}
\in\mathbb{Q}_{>0},
\]
and
\[
\deg_x O_{\ell-j(m-(\beta-1))}(x)<j(\alpha+1).
\]
Similarly, substituting the values of $f_{\ell-(j-1)(m-(\beta-1))}(x)$ and
$f_{\ell-(j-1)(m-(\beta-1))-\beta}(x)$ from Equations (\ref{eq3.6}) and (\ref{eq3.7}) in \Cref{eq3.5}, we obtain
\begin{align}
f'_{\ell-j(m-(\beta-1))-\beta}(x) &=(\ell-(j-1)(m-(\beta-1))-\beta)x^\alpha
\big[A_{\ell-(j-1)(m-(\beta-1))-\beta}\;x^{(j-1)(\alpha+1)}\notag \\
& +B_{\ell-(j-1)(m-(\beta-1))-\beta}\;x^{(j-2)(\alpha+1)+1}+O_{\ell-(j-1)(m-(\beta-1))-\beta}(x)\big]\notag \\
&-(\ell-(j-1)(m-(\beta-1)))\big[A_{\ell-(j-1)(m-(\beta-1))}x^{(j-1)(\alpha+1)} \notag\\
&+O_{\ell-(j-1)(m-(\beta-1))}(x)\big].\notag
\end{align}
Now integrating with respect to $x$, we get
\begin{align}
&f_{\ell-j(m-(\beta-1))-\beta}(x) = \frac{(\ell-(j-1)(m-(\beta-1))-\beta)A_{\ell-(j-1)(m-(\beta-1))-\beta}}{j(\alpha+1)}
x^{j(\alpha+1)}\notag \\
&\qquad+\frac{1}{(j-1)(\alpha+1)+1}\big[(\ell-(j-1)(m-(\beta-1))-\beta)B_{\ell-(j-1)(m-(\beta-1))-\beta}\notag \\
&\qquad-(\ell-(j-1)(m-(\beta-1)))A_{\ell-(j-1)(m-(\beta-1))}\big]x^{(j-1)(\alpha+1)+1}\notag\\
&\qquad+O_{\ell-j(m-(\beta-1))-\beta}(x),\notag
\end{align}
for some $O_{\ell-j(m-(\beta-1))-\beta}(x)\in k[x]$. \\ \noindent
By induction hypothesis we have $A_{\ell-(j-1)(m-(\beta-1))}\in \mathbb{Q}_{>0}, A_{\ell-(j-1)(m-(\beta-1))-\beta}\in k, B_{\ell-(j-1)(m-(\beta-1))-\beta}\in \mathbb{Q}_{<0}$,
$\deg_x O_{\ell-(j-1)(m-(\beta-1))}(x)<(j-1)(\alpha+1)$,
and $\deg_x O_{\ell-(j-1)(m-(\beta-1))-\beta}(x)<(j-2)(\alpha+1)+1$.\\
\noindent Hence it follows that 
\begin{align*}
f_{\ell-j(m-(\beta-1))-\beta}(x)
&=A_{\ell-j(m-(\beta-1))-\beta}x^{j(\alpha+1)}
+B_{\ell-j(m-(\beta-1))-\beta}x^{(j-1)(\alpha+1)+1}
\\
&+O_{\ell-j(m-(\beta-1))-\beta}(x),
\end{align*}
where
\[
A_{\ell-j(m-(\beta-1))-\beta}
=\frac{(\ell-(j-1)(m-(\beta-1))-\beta)A_{\ell-(j-1)(m-(\beta-1))-\beta}}{j(\alpha+1)}\in k,
\]
and
\begin{align}
&B_{\ell-j(m-(\beta-1))-\beta}=\big[
(\ell-(j-1)(m-(\beta-1))-\beta)B_{\ell-(j-1)(m-(\beta-1))-\beta}\notag \\
&-(\ell-(j-1)(m-(\beta-1)))A_{\ell-(j-1)(m-(\beta-1))}\bigg]\frac{1}{(j-1)(\alpha+1)+1}
\in\mathbb{Q}_{<0},\notag
\end{align}
with $\deg_x O_{\ell-j(m-(\beta-1))-\beta}(x)<(j-1)(\alpha+1)+1$.\\ 
Thus the claim follows.\\ 
\textbf{(II)} Finally, we show a contradiction from Lemma \ref{fval} and the previous claim.\\
Note that $j_0>1$. Substituting $j=j_0-1$ in Equations (\ref{eq3.6}) and (\ref{eq3.7}), we will get 
\begin{equation}\label{eq3.8}
f_{\ell-(j_0-1)(m-(\beta-1))}(x)=f_{m-\beta+1}
= A_{m-\beta+1}x^{(j_0-1)(\alpha+1)}+O_{m-\beta+1}(x).
\end{equation}
and
\begin{align}\label{eq3.9}
f_{\ell-(j_0-1)(m-(\beta-1))-\beta}(x)=f_{m-2\beta+1}
&= A_{m-2\beta+1}x^{(j_0-1)(\alpha+1)}
+ B_{m-2\beta+1}x^{(j_0-2)(\alpha+1)+1}
\notag\\
&+ O_{m-2\beta+1}(x).
\end{align}
where $A_{m-\beta+1}\in\mathbb{Q}_{>0}$, $A_{m-2\beta+1}\in k$, $B_{m-2\beta+1}\in\mathbb{Q}_{<0}$,
$O_{m-\beta+1}(x)\in k[x]$ with $\deg_x O_{m-\beta+1}(x)<(j_0-1)(\alpha+1)$, and $O_{m-2\beta+1}(x)\in k[x]$ with $\deg_x O_{m-2\beta+1}(x)<(j_0-2)(\alpha+1)+1$.

\medskip

\noindent Note that $m=(j+1)\beta$, for some $j\in \mathbb{N}$. Thus $m-\beta+1=j\beta+1$ and $m-2\beta+1=(j-1)\beta+1$ for $j\in \mathbb{N}$. \\
Suppose $g(x)=a_tx^t+a_{t-1}x^{t-1}+\dots+a_0$, where $a_i\in k\; \forall\; 0\leq i \leq t$ and $a_t\ne0$, be a polynomial of degree $t\geq0$.  
From Lemma \ref{fval}, we have
\begin{equation}\label{eq3.10}
f_{m-\beta+1}(x)=
\frac{x^{j\alpha} g(x)}{m-\beta+1},
\end{equation}
and
\begin{equation}\label{eq3.11}
f_{m-2\beta+1}(x)=
\frac{x^{(j-1)\alpha}g(x)}{m-2\beta+1}.
\end{equation}
Now, from Equations (\ref{eq3.8}) and (\ref{eq3.10}), we have 
$j\alpha+t=(j_0-1)(\alpha+1)\;\text{and}\; a_t=(m-\beta+1)A_{m-\beta+1}\in\mathbb{Q}_{>0}$. From Equations (\ref{eq3.9}) and (\ref{eq3.11}), we observe that, if  $A_{m-2\beta+1}\ne0$, then $(j-1)\alpha+t=(j_0-1)(\alpha+1)$, which contradicts the previous equality. Hence $A_{m-2\beta+1}=0$. Now, comparing the leading coefficient of $f_{m-2\beta+1}(x)$ from Equations (\ref{eq3.9}) and (\ref{eq3.11}), we get $a_t=(m-2\beta+1)B_{m-2\beta+1}\in\mathbb{Q}_{<0}$, which is again a contradiction.\\
Hence we get contradiction to the fact that $\ell>0$, and therefore $f\in k $, which implies that $g(x)=0$.
\end{proof}

\begin{prop}\label{mainprop}
    Let $k$ be a field of characteristic zero, and $\alpha\ge 1,\beta\ge 1$ and $m\ge \beta+1$ be integers. Consider the $k$-derivation $d=(1-y^\beta x^\alpha)\partial_{y}+y^m\partial_{x}$ of $k[x,y]$. If there exists some element $f\in k[x,y]$ such that $d(f)=af+g(x)$ for some  $a\in k^\ast $ and $g(x)\in k[x]$ with $\deg_x g(x)\le 1$, then $f\in k$ and $g(x)\in k$.
\end{prop}
\begin{proof}
       Let $a\in k^\ast$, and $g(x)\in k[x]$ with $\deg_x g(x)\le 1.$ Assume that there exists some element $f\in R_2\setminus k$ such that $d(f)=af+g(x).$ As $f\in R_2\setminus k$, we can write $f=\sum_{i=0}^\ell f_i(x)y^i$, where $\ell\geq 0$, $f_i(x)\in k[x]$ for all $0\leq i\leq \ell$ and $f_\ell(x)\neq 0$.\\
        Suppose $\ell=0$, that is, $f=f_0(x)\in k[x]$. Then $d(f)=af+g(x)$ implies that $y^mf'_0(x)=af_0(x)+g(x)$. Thus, $f'_0(x)=0.$ Therefore $f=f_0(x)\in k$, which is a contradiction. Hence $\ell>0$.\\
        Now, we are in the same setting as in Lemma \ref{deglem}, So we can use results from Lemma \ref{deglem}.\\
        From Lemma \ref{deglem}, we have 
        \begin{equation}\label{eq4:1}
            \deg_{x}f_{\ell-(m-(\beta-1))j_0}(x)=\deg_{x}f_0(x)=j_0(\alpha+1)
        \end{equation}
        and
        \begin{equation}\label{eq4:2}
            \deg_{x}f_{\ell-(m-(\beta-1))(j_0-1)}(x)=\deg_{x}f_{m-(\beta-1)}(x)=(j_0-1)(\alpha+1).
        \end{equation}
        Now, comparing the zero degree term in $y$ in \Cref{eq1:2}, we get
        $$f'_{-m}(x)+f_1(x)+(\beta-1)f_{-(\beta-1)}(x)-af_0(x)=g(x),$$
        \begin{equation}\label{eq4:3}
            \implies f_1(x)=af_0(x)+g(x).
        \end{equation}
        From \Cref{eq4:1}, we have $\deg_x f_0(x)=j_0(\alpha+1)$. Since  $j_0, \alpha\ge 1$, then $j_0(\alpha+1)\geq 2 >  \deg_x g(x)$. Thus, from \Cref{eq4:3}, it follows that $\deg_x f_1(x)=j_0(\alpha+1)$.\\

    \noindent \textbf{Claim:} We claim that for $1\le j\le m-\beta+1$, $$\deg_x f_j(x)=\left\lfloor \frac{j-1}{\beta} \right\rfloor \alpha+\deg f_1(x) \ \ \text{and} \ \ \delta_j=a ^{u-1} \delta_1 q_j,$$ where $q_j\in \mathbb{Q}_{>0},$ $j \equiv u \pmod{\beta}$ with $1\le u \le \beta$ and $\delta_s$ denotes the leading coefficient of $f_s(x)$ for all $1\leq s\leq m-\beta+1$.\\
    
    \noindent We prove the above claim by induction on $j$.\\
    The claim is trivially true for $j=1$. Let $t$ be a positive integer such that $1\leq t \leq m-\beta.$ Next, as the induction hypothesis, we assume for all integers $1\leq j\leq t$ that
    
    \textbf{(I)} $\deg_x f_j(x)=\left\lfloor\frac{j-1}{\beta}\right\rfloor\alpha+\deg_x f_1$ and 
    
    \textbf{(II)} $\delta_j=a^{u-1}\delta_1q_j$, where $j\equiv u\pmod{\beta}$ with $1\le u \le \beta$, $q_j\in\mathbb{Q}_{>0}$ and $\delta_s$ denotes the leading coefficient of $f
_s(x)$. \\

    \noindent Now we need to show that the claim holds for $t+1$.\\
    Substituting $i=t$ in \Cref{eq1:3}, we have that
    \begin{equation}\label{eq4:4}
        (t+1)f_{t+1}(x)=(t-\beta+1)x^\alpha f_{t-\beta+1}(x) +a f_t(x).
    \end{equation}
    Note that for $t\in \{1,\ldots, \beta-1\}$, the term $(t-\beta+1)x^\alpha f_{t-\beta+1}(x)=0$, and hence \Cref{eq4:4} becomes
    \begin{equation}\label{eq4:5}
        (t+1)f_{t+1}(x)=af_t(x),
    \end{equation}
      for all $1\le t \le \beta-1.$
    \noindent From the induction hypothesis, we have that $\deg_x f_t(x)=\left\lfloor \frac{t-1}{\beta}\right\rfloor \alpha+\deg_x f_1(x)$ and $\delta_t=a^{t-1}\delta_1 q_t$ for some $q_t\in \mathbb{Q}_{>0}$. Then, it follows from \Cref{eq4:5} that $\deg_x f_{t+1}(x)=\left\lfloor \frac{t-1}{\beta}\right\rfloor \alpha+\deg_x f_1(x)=\left\lfloor \frac{t}{\beta}\right\rfloor \alpha+\deg_x f_1(x)$ and $\delta_{t+1}=a^{t} \delta_1 q_{t+1}$, where $q_{t+1}=\frac{q_t}{t+1}\in \mathbb{Q}_{>0}.$ Thus, the claim holds for all $1\le t \le \beta-1.$ \\
    Henceforth, we assume that $t\geq \beta$.\\
    By (I) of the induction hypothesis, we have that $\deg_x f_{t}(x)=\left\lfloor \frac{t-1}{\beta} \right\rfloor \alpha+\deg f_1(x)$ and $\deg_x f_{t-\beta+1}(x)=\left\lfloor \frac{t-\beta}{\beta} \right\rfloor \alpha+\deg f_1(x)$. Then $\deg_x \ x^\alpha f_{t-\beta+1}(x)=\bigg(\left\lfloor \frac{t-\beta}{\beta} \right\rfloor +1\bigg) \alpha+\deg f_1(x)$ $= \left\lfloor \frac{t}{\beta} \right\rfloor \alpha+\deg f_1(x)$. Note that $\left\lfloor \frac{t-1}{\beta}\right\rfloor\leq \left\lfloor \frac{t}{\beta}\right\rfloor,$ for $t\geq \beta$. Therefore, we have that $$ \deg_x \ f_t(x)\leq \deg_x \ x^\alpha f_{t-\beta+1}(x).$$
    Now, we consider two distinct cases on $t$.\\
    \noindent \textbf{Case-1:} When $\beta 
    |t$, that is, $t=\beta v$ for some positive integer $v$. Note that in this case $\left\lfloor \frac{t-1}{\beta}\right\rfloor<\left\lfloor \frac{t}{\beta}\right\rfloor$. Therefore, $\deg_x \ f_t(x)< \deg_x \ x^\alpha f_{t-\beta+1}(x)$. Hence, from \Cref{eq4:4} it follows that, $\deg_x f_{t+1}=\left\lfloor \frac{t}{\beta}\right\rfloor\alpha +\deg_xf_1(x)$. By (II) of hypothesis, we have that $\delta_{t-\beta+1}=\delta_1 q_{t-\beta+1}$, where $q_{t-\beta+1}\in \mathbb{Q}_{>0}$. Thus from \Cref{eq4:4}, we have $\delta_{t+1}=\delta_1 q_{t+1}$, where $q_{t+1}=q_{t-\beta+1}\frac{t-\beta+1}{t+1}\in \mathbb{Q}_{>0}$.\\
    \textbf{Case-2:} When $t=\beta v +p$ where $v\in \mathbb{N}$ and $1\le p \le \beta-1$. Note that in this case $\left\lfloor \frac{t-1}{\beta}\right\rfloor= \left\lfloor \frac{t}{\beta}\right\rfloor$. Therefore, $\deg_x \ f_t(x)= \deg_x \ x^\alpha f_{t-\beta+1}(x)$. By (II) of hypothesis, we have that $\delta_{t-\beta+1}=a^p \delta_1 q_{t-\beta+1}$ and $\delta_{t}=a^{p-1}\delta_1 q_t$ where $q_{t-\beta+1},q_t\in \mathbb{Q}_{>0}$. Then, from \Cref{eq4:4}, we have that
    $$ (t+1)\delta_{t+1}=a^p (t-\beta+1)\delta_1 q_{t-\beta+1}+a\cdot a^{p-1}\delta_1 q_t=a^p\delta_1((t-\beta+1)q_{t-\beta+1}+q_t).$$
    Since $q_t, q_{t-\beta+1},t-\beta+1\in \mathbb{Q}_{>0}$, then $(t-\beta+1)q_{t-\beta+1}+q_t\in \mathbb{Q}_{>0}$. Hence $\delta_{t+1}=a^p \delta_1 q_{t+1}$, where $q_{t+1}=\frac{(t-\beta+1)q_{t-\beta+1}+q_t}{t+1}\in \mathbb{Q}_{>0}$ and $\deg_x f_{t+1}=\left\lfloor \frac{t}{\beta}\right\rfloor \alpha +\deg_x f_1(x)$. 
    This proves the claim.\\
    In particular, we have $\deg_x f_{m-\beta+1}(x)=\left\lfloor \frac{m-\beta}{\beta} \right\rfloor \alpha+\deg_x f_1(x)$. Recall that $\deg_x f_1(x)$ $= j_0(\alpha+1).$ Thus, $\deg_x f_{m-\beta+1}(x)=\left\lfloor \frac{m-\beta}{\beta} \right\rfloor \alpha+ j_0(\alpha+1)$. However, we have from \Cref{eq4:2} that $\deg_x f_{m-\beta+1}(x)=(j_0-1)(\alpha+1)$. Therefore, 
    \begin{align*}
        \left\lfloor \frac{m-\beta}{\beta} \right\rfloor \alpha+ j_0(\alpha+1)&=(j_0-1)(\alpha+1),\\
         \implies \left\lfloor\frac{m-\beta}{\beta}\right\rfloor\alpha&=-(\alpha+1),
    \end{align*}
    which is absurd as $m\ge \beta +1$ and $\alpha\geq 1$.
    Hence, we get a contradiction to the fact that $\ell>0$. Thus $f\in k$, and therefore $d(f)=af+g(x)=0$, which implies that $g(x)=-af\in k$.
\end{proof}

\section{The General Case}\label{The General Case}

	Let $k$ be a field of characteristic zero, and let $n>2$ be a positive integer. Set $R_2 = k[x_1,x_2]$ and $R_n = R_2[x_3, \ldots, x_n]$. 
Suppose $d_2$ is a $k$-derivation of $R_2$, and define a $k$-derivation $d_n$ of $R_n$ by

$$d_n = d_2 + (a_2 x_{3} + g_2(x_2))\partial_{x_{3}}+\cdots + (a_{n-1} x_n + g_{n-1}(x_{n-1}))\partial_{x_n},$$

where $a_i \in k$ and  $g_i(x_i) \in k[x_i] \setminus k$ for all $2 \leq i \le n-1$. Note that $d_n|_{R_2}=d_2$ and $d_2(x_i)=0$ for $3\leq i\leq n$. Let $j$ be any integer such that $2 <j \leq n$, and define $R_j = R_2[x_3, \ldots, x_j]$, by $d_j$ we denote the restriction of $d_n$ to $R_j$; explicitly, $d_j:=d_2 + (a_2 x_3 + g_2(x_2))\partial_{x_3} + \cdots + (a_{j-1} x_j + g_{j-1}(x_{j-1}))\partial_{x_j}$. Thus, $d_j$ is a $k$-derivation of $R_j.$

\begin{lem}\label{lemma no units}
  Let $n$, $R_2$, $d_2$ and $d_n$ be as above. Further, assume that $\deg_{x_2} g_2(x_2) \leq l$ for some $l \in \mathbb{N} \cup \{\infty\}$, and that $d_2$ satisfies the following condition:
	
	If $d_2(r) = a r + g(x_2)$ for some $r \in R_2$, $a \in k$, and $g(x_2) \in k[x_2]$ with $\deg_{x_2} g(x_2) \leq l$, 
	then $r \in k$ and $g \in k$.
    
 Let $j$ be any integer such that $2 <j \leq n$. If $d_n(f) = c f + h(x_j)$ for some $f\in R_j$, $c \in k$ and $h(x_j) \in k[x_j]$, then $f \in k$ and $h \in k$.
\end{lem}

\begin{proof} Note that for any $f\in R_j$, for $3\leq j\leq n$, we have $d_n(f)=d_j(f)$.
	
	(i) We first prove the result for $j=3$. Recall that $d_{3}=d_2+(a_2x_{3}+g_2(x_2))\partial_{x_3}$ where $a_2\in k$, $g_2(x_2)\in k[x_2]\setminus k$ with $\deg_{x_2}g_2(x_2)\leq l$. Let $f\in R_{3}$ such that $d_{3}(f)=cf+h(x_{3})$ for some $c \in k$ and $h(x_{3})\in k[x_{3}]$.  If $f=0$, we are done. Otherwise $f\in R_{3}\backslash\{0\}$, and we can write $f=\sum_{j=0}^{t}f_jx_{3}^j$ where $t\geq 0$ and $f_j\in R_2$ for all $0\leq j\leq t$ and $f_t\neq 0$.
	If $t=0$, that is $f\in R_{2}$, then $d_{3}(f)=d_{2}(f)$. Now, $d_{3}(f)=d_{2}(f)=cf+h(x_{3})$ and $d_2(f)\in R_2$ implies that $h(x_{3})\in k$. Since $d_{2}(f)=cf+h$, using the property of $d_2$, we get that $f\in k$. Thus, for $t=0$ the result holds.
    
	Now, we suppose $t\geq 1$, Then
	\begin{align}\label{onlyeqn}
		d_{3}(f)&=d_{3}(\sum_{j=0}^{t}f_jx_{3}^j)\\ 
		&=(d_{3}(f_t)+a_2tf_t)x_{3}^t+(d_{3}(f_{t-1})+g_2(x_2)tf_t+a_2(t-1)f_{t-1})x_{3}^{t-1}+ T, \nonumber
	\end{align} 
	where $a_2\in k$, $g_2(x_2)\in k[x_2]\setminus k$, $f_t\neq 0$ and $T \in R_{2}[x_{3}]$ with $\deg_{x_{3}}(T) <t-1$.
	
	Now, write $d_{3}(f)=cf+h(x_{3})$ as \begin{equation}\label{eqn-d_{3}}
	    d_{3}(f)-cf=h(x_{3}).
	\end{equation} Then, simplify \Cref{eqn-d_{3}} by substituting $f=\sum_{j=0}^{t}f_jx_{3}^j$ and the expression of $d_{3}(f)$ from \Cref{onlyeqn}, we get that
    \begin{align}\label{eqn-two}
        (d_{3}(f_t)+a_2tf_t-cf_t)x_{3}^t+(d_{3}(f_{t-1})+g_2(x_2)tf_t+a_2(t-1)f_{t-1}-cf_{t-1})x_{3}^{t-1} \nonumber\\
        + T=h(x_{3}).
    \end{align}
    Since $h(x_{3})\in k[x_{3}]$, on comparing the coefficients of $x_{3}^t$ and $x_{3}^{t-1}$ in \Cref{eqn-two},we get \begin{equation}\label{one}
		d_{3}(f_t)+a_2tf_t-cf_t\in k \end{equation} and
	\begin{equation}\label{two}
 d_{3}(f_{t-1})+g_2(x_2)tf_t+a_2(t-1)f_{t-1}-cf_{t-1}\in k
	\end{equation} respectively.
	Note that for all $j$, $ 0\leq j\leq t$, we have $d_{3}(f_{j})=d_{2}(f_{j})$ since $f_j\in R_{2}$. In particular $d_2(f_t)=d_{3}(f_t)$ and $d_{2}(f_{t-1})=d_{3}(f_{t-1})$. From \Cref{one}, we get that 
	\begin{equation}\label{eqn lem no units}
	    d_2(f_t)=d_{3}(f_t)=(c-a_2t)f_t+\alpha,
	\end{equation}
 for some $\alpha \in k$ and then by the assumption on $d_2$, we get that $f_t\in k^{\ast}$. Also, using $f_t\in k^{\ast}$ in \Cref{two} we get that \[d_{2}(f_{t-1})=d_{3}(f_{t-1})=(c-a_2(t-1))f_{t-1}-g_2(x_2)tf_t+\beta,\] for some $\beta \in k$. 
	Since $c-a_2(t-1)\in k$ and $-g_2(x_2)tf_t+\beta\in k[x_2]$ with $\deg_{x_2}(-g_2(x_2)tf_t+\beta)\leq l$, by the assumption on $d_2$, it follows that $f_{t-1}\in k$ and thus $d_{3}(f_{t-1})=0$. Then $g_2(x_2)=\beta(tf_t)^{-1}+(c-a_2(t-1))((tf_t)^{-1}f_{t-1})\in k$, which is a contradiction as $g_2(x_2)\in k[x_2]\setminus k$. 
	
	(ii) Note that it is enough to prove the result for $j=4$. From (i), we know that if $d_3(f)=cf+h(x_3)$ for some $f\in R_3$, $c\in k$ and $h(x_3)\in k[x_3]$ then $f\in k$ and $h(x_3)\in k.$ Now consider $d_4=d_3+(a_3x_4+g_3(x_3))\partial_{x_4}.$ Applying the same argument as in (i), with $d_3,a_3,g_3(x_3)$ and $x_4$ in place of $d_2,a_2,g_2(x_2)$ and $x_3$, respectively, we get that the Lemma holds for $j=4.$
\end{proof}  

\begin{rem}
    Let $a=0$ in the hypothesis of Lemma~\ref{lemma no units}. We have from the proof of the Lemma~\ref{lemma no units}, the \Cref{eqn lem no units}, that is,  $$d_2(f_t)=d_{3}(f_t)=(c-a_2t)f_t+\alpha.$$ If $c=0,$ that is, $d_n(f)=h(x_j)$ for $f\in R_j$ and $h(x_j)\in k[x_j]$ then in order to apply the hypothesis on $d_2$, in Equation~\ref{eqn lem no units}, it is necessary that $a_2=0.$ Similarly, we have the same observation in each subsequent extension step, and hence we get $a_i=0$ for all $2\leq i\leq n-1$. Thus, we conclude that when $c=0$, all $a_i$ must be zero to apply the hypothesis with $a=0.$
\end{rem}

Let $n,a_i$ and $g_i(x_i)$ be as in Lemma~\ref{lemma no units}. We consider two classes of $k$-derivations of $R_n=k[x_1,\ldots,x_n]$, namely, the $k$-derivation $$(1 +x_1^{m_2}x_2)\partial_{x_1}+x_1^{m_1}\partial_{x_2}+(a_2 x_3 +g_2(x_2))\partial_{x_3}+\cdots+(a_{n-1}x_n + g_{n-1}(x_{n-1}))\partial_{x_n},$$ where $m_1$ and $m_2$ be two positive integers such that $m_2$ does not divide $m_1$, and the $k$-derivation $$(1-x_2^{\alpha} x_1^{\beta})\partial_{x_1} +x_1^m\partial_{x_2}+(a_2 x_3 +g_2(x_2))\partial_{x_3}+\cdots+(a_{n-1}x_n + g_{n-1}(x_{n-1}))\partial_{x_n},$$ where $\alpha\geq 1,\beta\geq 1$ and $m\geq \beta+1$ be integers. 

For $n=2,$ letting $x_1=y$ and $x_2=x$, we get the ring $R_2=k[x,y]$ with the $k$-derivations $y^{m_1}\partial_{x}+(1+y^{m_2}x)\partial_{y}$ and $y^m\partial_x+(1-x^{\alpha} y^{\beta})\partial_y$, which have been considered in the previous section.
\begin{cor}
    Let $m_1$ and $m_2$ be two positive integers such that $m_2$ does not divide $m_1$. Let $d_2=(1 +x_1^{m_2}x_2)\partial_{x_1}+x_1^{m_1}\partial_{x_2}.$ Let $n\geq 3.$ Let $a_i$ and $g_i(x_i)$, for all $2\leq i\leq n-1$, be as in Lemma~\ref{lemma no units}, and let $d_n=(1 +x_1^{m_2}x_2)\partial_{x_1}+x_1^{m_1}\partial_{x_2}+(a_2 x_3 +g_2(x_2))\partial_{x_3}+\cdots+(a_{n-1}x_n + g_{n-1}(x_{n-1}))\partial_{x_n}$. Let $j$ be any integer such that $2 <j \le n$. Now, if $d_n(f) = c f + h(x_j)$ for some $f\in R_j$, $c \in k$ and $h(x_j) \in k[x_j]$, it follows that $f \in k$ and $h \in k$. 
\end{cor}
\begin{proof}
    By Proposition~\ref{prop:main}, $d_2$ satisfies the hypothesis of Lemma~\ref{lemma no units} for $l=\infty$. Then the proof follows using Lemma~\ref{lemma no units}.
\end{proof}

\begin{cor}
    Let $\alpha,\beta\geq 1$ and $m\geq \beta+1$ be integers. Let $d_2=(1-x_2^{\alpha} x_1^{\beta})\partial_{x_1} +x_1^m\partial_{x_2}.$ Let $n\geq 3.$ Let $a_i$ and $g_i(x_i)$, for all $2\leq i\leq n-1$, be as in Lemma~\ref{lemma no units}, and let $d_n=(1-x_2^{\alpha} x_1^{\beta})\partial_{x_1} +x_1^m\partial_{x_2}+(a_2 x_3 +g_2(x_2))\partial_{x_3}+\cdots+(a_{n-1}x_n + g_{n-1}(x_{n-1}))\partial_{x_n}$. Let $j$ be any integer such that $2 <j \leq n$. Assume that $d_n(f) = c f + h(x_j)$ for some $f\in R_j$, $c \in k$ and $h(x_j) \in k[x_j]$. Then, 
    \begin{enumerate}
        \item[(i)] if $c=0$ and $a_i=0$ for all $2\leq i\leq j-1$  then $f \in k$ and $h(x_j) =0$, and  
        \item[(ii)] if $\deg_{x_2}g_2(x_2)\leq 1$ then $f \in k$ and $h \in k.$
    \end{enumerate}
\end{cor}
\begin{proof}
  (i) Let $c=0.$ Then $d_n(f)= h(x_j)$ for some $f\in R_j$ and $h(x_j) \in k[x_j]$. By Proposition~\ref{pc=0}, we have that if $d_2(r)=g(x_2)$ for some $r\in R_2$ and $g(x_2)\in k[x_2]$ then $r\in k$ and  $g(x_2)=0$. Thus, the hypothesis of the Lemma~\ref{lemma no units} is satisfied with $a=0$ and $l=\infty.$ Then it follows using  Lemma~\ref{lemma no units} that $f\in k$ and $h(x_j)\in k.$ Since $f\in k,$ $d_n(f)=0.$ Hence $h(x_j)=0,$ proving (i).

  (ii) Let $\deg_{x_2}g_2(x_2)\leq 1$. By Proposition~\ref{pc=0} and Proposition~\ref{mainprop}, if $d_2(r)=ar+g(x_2),$ where $a\in k$ and $\deg_{x_2}g(x_2)\leq 1$ then $r\in k$ and $g(x_2)\in k.$ Therefore, $d_2$ satisfies the hypothesis of Lemma~\ref{lemma no units} with $l=1$. Then it follows from Lemma~\ref{lemma no units} that $f\in k$ and $h(x_j)\in k,$ proving (ii).
\end{proof}

Next, we recall a result by Shamsuddin from \cite{Sham77}.
\begin{thm}[Shamsuddin' result]\label{Sham result}
 Let $R$ be a ring containing $\mathbb{Q}$ and let $d$ be a simple derivation of $R.$ Extend the derivation $d$ to a derivation $\tilde{d}$ of the polynomial ring $R[t]$ by setting $\tilde{d}(t)=at+b$ where $a,b\in R.$ Then the following conditions are equivalent:
 \begin{enumerate}
     \item[(i)] $\tilde{d}$ is simple.
     \item[(ii)] There exist no elements $r$ of $R$ such that $d(r)=ar+b.$
 \end{enumerate}
\end{thm}
\begin{cor}\label{cor simple}
    Let $n,R_2,R_n,d_2$ and $d_n$ be as defined above. Let $\deg_{x_2}g_2(x_2)\leq l$ for some $l \in \mathbb{N} \cup  \{\infty\}.$ Suppose that $d_2$ satisfies the following condition: If $d_2(r)=ar+g(x_2)$ for some $a\in k,$ $r\in R_2$ and $g(x_2)\in k[x_2]$ with $\deg_{x_2}g(x_2)\leq l,$ then $r\in k$ and $g(x_2)\in k.$ If $d_2$ is a simple derivation of $R_2$ then $d_n$ is a simple derivation of $R_n.$
\end{cor}
\begin{proof}
    It is enough to prove the result for $d_3$, then the rest of the proof follows similarly using Lemma~\ref{lemma no units}. Recall that $d_3=d_2+(a_2x_3+g_2(x_2))\partial_{x_3},$ where $g_2(x_2)\in k[x_2]\setminus k$ and $\deg_{x_2}g_2(x_2)\leq l.$ Thus $d_3(x_3)=a_2x_3+g_2(x_2)$. Now, let if possible, there exists an $r\in R_2$ satisfying $d_2(r)=a_2r+g_2(x_2)$. Since $\deg_{x_2}g_2(x_2)\leq l,$ it follows from hypothesis that $r\in k$ and $g_2(x_2)\in k.$ But $g_2(x_2)\notin k$ and hence we get a contradiction. Then the derivation $d_3$ is simple by Theorem~\ref{Sham result}. 
\end{proof}
\begin{rem}
    Let $a=0$ in the hypothesis of Corollary~\ref{cor simple}. Then, in the proof of Corollary~\ref{cor simple}, to show that $d_3$ is simple, that is, to apply the hypothesis to $d_2$ with $a=0$, in the equation $d_2(r)=a_2r+g_2(x_2),$ we must have $a_2=0$. Similarly, to show that $d_4$ is simple, we require $a_3=0$, and so on. Thus, Corollary~\ref{cor simple} holds true with the restriction $a=0$ in the hypothesis provided $a_i=0$ for all $2\leq i \leq n-1.$  
\end{rem}
\begin{lem}\label{lem simple}
    Let $m,m_1,m_2,\alpha$, and $\beta$ be positive integers. Then $k$-derivations $(1+x_1^{m_2}x_2)\partial_{x_1}+x_1^{m_1}\partial_{x_2}$ and $(1-x_2^{\alpha} x_1^{\beta})\partial_{x_1} +x_1^m\partial_{x_2}$ are simple derivations of $k[x,y].$
\end{lem}
\begin{proof}
    It follows from \cite[Theorem 2.6]{SP25}. 
\end{proof}
\begin{rem}
The simplicity of derivations in Lemma~\ref{lem simple} can also be derived from \cite[Theorem 4.1 and Theorem 6.1]{Kour14}.
\end{rem}
\begin{cor}
     Let $m_1$ and $m_2$ be two positive integers such that $m_2$ does not divide $m_1$. Let $d_2=(1 +x_1^{m_2}x_2)\partial_{x_1}+x_1^{m_1}\partial_{x_2}.$ Let $n\geq 3.$ Let $a_i$ and $g_i(x_i)$, for all $2\leq i\leq n-1$, be as in Lemma~\ref{lemma no units}, and let $d_n=(1 +x_1^{m_2}x_2)\partial_{x_1}+x_1^{m_1}\partial_{x_2}+(a_2 x_3 +g_2(x_2))\partial_{x_3}+\cdots+(a_{n-1}x_n + g_{n-1}(x_{n-1}))\partial_{x_n}$. Then $d_n$ is simple.
\end{cor}
\begin{proof}
    By Lemma~\ref{lem simple}, $d_2$ is simple. By Proposition~\ref{prop:main}, $d_2$ satisfies the hypothesis of Corollary~\ref{cor simple}. The proof then follows from Corollary~\ref{cor simple}. 
\end{proof}
\begin{cor}
     Let $\alpha,\beta\geq 1$ and $m\geq \beta+1$ be integers. Let $d_2=(1-x_2^{\alpha} x_1^{\beta})\partial_{x_1} +x_1^m\partial_{x_2}.$ Let $n\geq 3.$ Let $a_i$ and $g_i(x_i)$, for all $2\leq i\leq n-1$, be as in Lemma~\ref{lemma no units}, and let $d_n=(1-x_2^{\alpha} x_1^{\beta})\partial_{x_1} +x_1^m\partial_{x_2}+(a_2 x_3 +g_2(x_2))\partial_{x_3}+\cdots+(a_{n-1}x_n + g_{n-1}(x_{n-1}))\partial_{x_n}$. Then $d_n$ is simple under each of the following assumptions:
     \begin{enumerate}
         \item[(i)] if $a_i=0$ for all $2\leq i \leq n-1$,
         \item[(ii)] if $\deg_{x_2}g_2(x_2)\leq 1.$
     \end{enumerate}
    
\end{cor}
\begin{proof}
     By Lemma~\ref{lem simple}, $d_2$ is simple. Consider the assumption $a_i=0$ for all $2\leq i \leq n-1$, we have from Proposition~\ref{pc=0} that $d_2$ satisfies the hypothesis of Corollary~\ref{cor simple} when $a=0, l=\infty$. The proof then follows from Corollary~\ref{cor simple}. Now assume that $\deg_{x_2}g_2(x_2)\leq 1,$ we have from Proposition~\ref{pc=0} and Proposition~\ref{mainprop} that $d_2$ satisfies the hypothesis of Corollary~\ref{cor simple} when $\deg_{x_2}g(x_2)\leq 1$, that is, $l=1$. The proof then follows from Corollary~\ref{cor simple}. 
\end{proof}

\end{document}